\documentclass[11pt]{article}
\usepackage[english]{babel}
\usepackage{amsmath,amssymb,amsfonts}
\newtheorem{theorem}{Theorem}[section]
\newtheorem{lemma}[theorem]{Lemma}
\newtheorem{proposition}[theorem]{Proposition}
\newtheorem{definition}{Definition}[section]

\newtheorem{hypothesis}[theorem]{Hypothesis}

\newtheorem{remark}[theorem]{Remark}

\def\sqr#1#2{{\vcenter{\vbox{\hrule height .#2pt
     \hbox{\vrule width .#2pt height#1pt \kern#1pt \vrule
     width .#2pt} \hrule height .#2pt}}}}
\def\square{\mathchoice\sqr54\sqr54\sqr{4.1}3\sqr{3.5}3}

\makeatletter
\@addtoreset{equation}{section}

\makeatother

\def\qed{{\hfill\hbox{\enspace${ \square}$}} \smallskip}
\def\sqr#1#2{{\vcenter{\vbox{\hrule height .#2pt \hbox{\vrule
 width .#2pt height#1pt \kern#1pt \vrule
width .#2pt} \hrule height .#2pt}}}}
\def\square{\mathchoice\sqr54\sqr54\sqr{4.1}3\sqr{3.5}3}

\def\ds{\begin{displaystyle}}
\def\eds{\end{displaystyle}}

\def\<{\langle }
\def\>{\rangle }

\title{Nonlinear Backward Stochastic Evolutionary Equations\\ Driven by a Space-Time White Noise}
\author{Ying Hu\thanks{%
IRMAR, Universit{\'e} Rennes 1, Campus de Beaulieu, 35042 Rennes
Cedex, France, and School of Mathematical Sciences, Fudan
University, Shanghai 200433, China. Partially supported by Lebesgue
Center of Mathematics ``Investissements d'avenir"
program-ANR-11-LABX-0020-01, by ANR CAESARS (Grant No. 15-CE05-0024) and by
ANR MFG (Grant No. 16-CE40-0015-01). email: \texttt{ying.hu@univ-rennes1.fr}
} \and
Shanjian Tang\thanks{%
Department of Finance and Control Sciences, School of Mathematical
Sciences, Fudan University, Shanghai 200433, China. Partially
supported by National Science Foundation of China (Grant No.
11631004) and Science and Technology Commission of Shanghai
Municipality (Grant No. 14XD1400400). email:
\texttt{sjtang@fudan.edu.cn}}}

\begin{document}
\maketitle

\begin{abstract}
We study the well solvability of nonlinear backward stochastic evolutionary equations driven by a space-time white noise. We first establish a novel a priori estimate for solution of linear backward stochastic evolutionary equations, and then give an existence and uniqueness result for nonlinear backward stochastic evolutionary equations. A dual argument plays a crucial role in the proof of these results. Finally, an example is given to illustrate the existence and uniqueness result.
\end{abstract}

\section{Introduction}
Let $H$ be a Hilbert space with $\{e_i\}$ being its  orthonormal basis,  $A$ an infinitesimal generator which generates a strongly continuous semigroup $\{e^{At}, t\ge0\}$, and  $\mathcal{S}_2(H)$
 the Hilbert space of Hilbert-Schmidt operators in $H$. Denote by $W$ a cylindrical Wiener process in the probability space $(\Omega, \mathcal{F}, \mathbb{P})$ with $(\mathcal{F}_t)_{t\geq 0}$ being
 the augmented natural filtration and $\mathcal{P}$ the predictable $\sigma$-algebra.

By $L^0_{\mathcal{P}}([0,T]\times \Omega,H)$ we denote the totality of $H$-valued progressively measurable processes $X$. For $p\in [1,\infty)$, by $L^p_{\mathcal{P}}([0,T]\times \Omega,H)$ we denote the Banach space of $H$-valued progressively measurable processes $X$ with $\mathbb{E}\int_0^T \|X_s\|^pds < \infty$ and by $L^p_{\mathcal{P}}(\Omega,{C}([0,T],H))$ the subspace  of $H$-valued progressively measurable processes $X$ with strongly continuous trajectories satisfying $\mathbb{E}\max_{s\in[0,T]}\|X_s\|^p < \infty$.
Here and below we use the symbol $\|\cdot\|$ to denote a norm when  the corresponding
space is clear from the context, otherwise we use a subscript. Consider the map $f: \Omega\times [0,T]\times H\times \mathcal{S}_2(H)\times H\to H$ where $f(\cdot,0,0,0)\in L^0_{{\cal P}}(\Omega\times [0,T],H)$ and there is a positive constant $L$ such that
$$
\|f(t,p_1,q_1,s_1)-f(t,p_2,q_2,s_2)\|_H\le L(\|p_1-p_2\|_H+\|q_1-q_2\|_{{\mathcal S}_2(H)}+\|s_1-s_2\|_H).
$$

Linear backward stochastic evolutionary equations arise in the formulation of stochastic maximum principle for optimal control of stochastic partial differential equations, and see among others ~\cite{Be,HuPe1,DuMe2,FuHuTe,Gu,LuZh,TaLi,Zh}. The study can be dated back to the work of Bensoussan~\cite{Be}, and  to Hu and Peng~\cite{HuPe1} for a general context. The nonlinear case is given by Hu and Peng~\cite{HuPe2}. In these works, the underlying Wiener process is assumed to have a trace-class covariance operator---in particular, to be  finite-dimensional. Recently, Fuhrman, Hu, and Tessitore~\cite{FHT} discusses a linear  backward stochastic evolutionary equation driven by a space-time white noise. The objective of the paper is to study the nonlinear backward stochastic evolutionary equation driven by a space-time white noise.

Consider the following form of nonlinear backward stochastic evolutionary equations (BSEEs)
\begin{equation}\label{nlbsde}
    \left\{\begin{array}{lll}
    -dP_t&=&\displaystyle \left[A^*P_t + \sum_{i=1}^\infty C_i^*(t)Q_te_i+f\left(t,P_t,Q_t,\sum_{i=1}^\infty C_i^*(t)Q_te_i\right) \right]\, dt\\
    &&\displaystyle -\sum_{i=1}^\infty Q_te_i\,d\beta_t^i, \quad t\in (0,T];   \\
    P_T&=&\eta,
\end{array}\right.
\end{equation}
where $\beta^i_t=\< e_i, W_t\>$, $i=1,2...$ is a family of independent Brownian motions,
  $\eta\in L^2(\Omega,{\cal F}_T,\mathbb P,H)$.
The unknown process is the pair denoted $(P,Q)$ and takes values
in $H\times \mathcal{S}_2(H)$.
We will work under the following assumptions, which are assumed to
hold throughout the paper.

\medskip

\begin{hypothesis}\label{Hyp:C_i} $ $
\begin{enumerate}

  \item  $e^{tA}$, $t\ge0$, is a strongly continuous semigroup of
  bounded linear operators in $H$. Moreover,
  $ e^{tA}\in {\mathcal{S}_2(H)}$ for all $ t> 0$ and there exist
  constants $c>0$ and $\alpha \in [0, 1/2) $ such that
  $\|e^{tA}\|_{\mathcal{S}_2(H)} \leq c t^{-\alpha}$ for all $ t\in ( 0,T]$.

\item The processes $C_i $ are strongly progressively
measurable with values in $\mathcal{L}(H)$. Moreover
we have $ \|C_i(t)\|_{\mathcal{L}(H)}\leq c,\quad \hbox{$\mathbb{P}$-a.s.}$ for all $t\in [0,T]$
and $i\in \mathbb{N}$.

  \item $\sum_{i=1}^{\infty} \|e^{tA}C_i(s) h\|^2 \leq c t^{-2\alpha}\|h\|_H^2$ for
  all $ t\in ( 0,T],$ $s\geq 0$, and $h\in H$.

\end{enumerate}
\end{hypothesis}

\medskip
We give the following notion of (mild) solution.

 \begin{definition}\label{def:solution_BSDE} We say that a pair of processes
 $$(P,Q)\in  L^2_{{\cal P}}(\Omega\times [0,T],H)\times L^2_{{\cal P}}(\Omega\times [0,T],\mathcal{S}_2(H))$$
  is a mild solution to equation \eqref{nlbsde} if the following holds:

  \begin{enumerate}

\item The sequence
 $$S^M(s) : =\sum_{i=1}^M (T-s)^{\alpha}  C_i^*(s)Q_se_i, \quad s\in [0,T]; \quad M=1,2,\ldots,$$
converges weakly in $L^{2}_{\mathcal{P}}(\Omega \times [0,T],H)$.

\item For any $t\in [0,T]$,
\begin{equation}\label{eq:mild}
\begin{array}{rcl}
P_t&=&\displaystyle e^{(T-t)A^*}\eta+ \sum_{i=1}^\infty \int_t^T e^{(s-t)A^*} C_i^*(s)Q_se_ids\\
&&\displaystyle+\int_t^T e^{(s-t)A^*}f\left(s,P_s,Q_s,\sum_{i=1}^\infty C_i^*(s)Q_se_i\right)\, ds\\
&&\displaystyle -\sum_{i=1}^\infty\int_t^T e^{(s-t)A^*}Q_s e_id\beta_s^i,\, \hbox{$\mathbb{P}$-a.s.}
\end{array}
\end{equation}
\end{enumerate}
 \end{definition}

Note that the integral involving $f$ is well-defined due to the following
\begin{eqnarray}
  &&\int_t^T \left\|e^{(s-t)A^*}f\left(s,P_s,Q_s,\sum_{i=1}^\infty C_i^*(s)Q_se_i\right)\right\|_H\, ds \nonumber\\
   &\le&  C \int_t^T \left(\|f(s,0,0,0)\|_H+\|P_s\|_H+\|Q_s\|_{{\mathcal S}_2(H)}+(T-s)^{-\alpha}\left\|(T-s)^\alpha\sum_{i=1}^\infty C_i^*(s)Q_se_i\right\|_H\right)\, ds\nonumber\\
   &<&\infty. \nonumber
\end{eqnarray}

Note that the term $\sum_{i=1}^\infty C_i^*Qe_i$ is not bounded in $Q$ with respect to the Hilbert-Schmidt norm. Its appearance in the drift gives rise to new difficulty in the resolution of the underlying BSEEs, and has to be carefully estimated. In particular, we prove via a dual method novel a priori estimate (see Proposition~\ref{apriori} in Section 2 below for details) for solution of linear BSEEs driven by a space-time white noise. The new a priori estimate and the dual arguments are  crucial in the subsequent Picard iteration for our nonlinear BSEEs.

The rest of the paper is organized as follows. In Section 2, we prove a new a priori estimate for solution of linear BSEEs driven by a space-time white noise. The nonlinear BSEEs are studied in Section 3. Finally in Section 4, we give an example.

\section{Linear BSEE Revisited}

(Forward) stochastic evolutionary equations (FSEE) driven by a cylindrical Wiener process have been extensively studied. See, e.g. Da Prato and Zabczyk~\cite{DaZa,DaZa2} for excellent expositions and the references therein. Here we give some precise a priori estimate for mild solutions of linear FSEEs, which will play a crucial role in our subsequent analysis.

\begin{lemma} \label{2.2}
For $\gamma\in L^\infty_{\cal P}(\Omega\times [0,T],H)$, the linear stochastic equation
\begin{equation}\label{lfsde}
    \left\{\begin{array}{lll}
    d\mathcal{Y}_t&=&\displaystyle A\mathcal{Y}_t\,dt + \sum_{i=1}^\infty C_i(t)\mathcal{Y}_t\,d\beta_t^i+\sum_{i=1}^\infty C_i(t)(T-t)^\alpha\gamma_t\,d\beta_t^i,\quad t\in (0,T];
    \\
    \mathcal{Y}_0&=&0,
\end{array}\right.
\end{equation}
has a unique mild solution $\mathcal{Y}^\gamma$ in $L^2_{\cal P}(\Omega, C([0,T],H)).$ Furthermore, we have
\begin{equation}\label{max estimate}
    \mathbb{E}\|\mathcal{Y}^\gamma_t\|^2_H\le C  \mathbb{E}\int_0^t(t-s)^{-2\alpha}(T-s)^{2\alpha}\|\gamma_s\|^2_H\, ds.
\end{equation}
For $\eta\in L^2(\Omega, \mathcal{F}_T,\mathbb{P}, H)$ and $f_0\in L^0_{\cal P}(\Omega\times [0,T],H)$ such that for some $\beta\in (0,\frac{1}{2})$
\begin{equation}\label{f0}
   {\mathbb E}\int_0^T(T-s)^{2\beta}\|f_0(s)\|^2_H\, ds<\infty,
\end{equation}
the following linear functional $G$ defined by
\begin{equation}\label{LF}
  G(\gamma):=E\langle \eta, \mathcal{Y}^\gamma_T\rangle+E\int_0^T\langle f_0(t), \mathcal{Y}^\gamma_t\rangle\, dt, \quad \gamma\in L^\infty_{\cal P}(\Omega\times [0,T],H),
\end{equation}
has a unique  linear and continuous extension to $L^2_{\cal P}(\Omega\times [0,T],H)$, denoted by  $\overline G$.
\end{lemma}

\noindent{\bf Proof. } The first two assertions can be found in \cite[Theorem 4.3 and Proposition 4.5]{FHT}. We now prove the last assertion, that is
\begin{eqnarray}
|G(\gamma)|^2&\le& C\|\gamma\|^2_{L^2_{\cal P}(\Omega\times [0,T],H)}, \quad \gamma\in L^\infty_{\cal P}(\Omega\times [0,T],H).
\end{eqnarray}

First from \eqref{max estimate}, we immediately have
$$\mathbb{E}\|\mathcal{Y}^\gamma_T\|^2_H\le C \|\gamma\|^2_{L^2_{\cal P}(\Omega\times [0,T],H)}, \quad \gamma\in L^\infty_{\cal P}(\Omega\times [0,T],H). $$
It suffices to prove the following
$$\mathbb{E}\int_0^T (T-s)^{-2\beta}\|\mathcal{Y}^\gamma_s\|^2_H\, ds\le C \|\gamma\|^2_{L^2_{\cal P}(\Omega\times [0,T],H)}, \quad \gamma\in L^\infty_{\cal P}(\Omega\times [0,T],H). $$

We have from \eqref{max estimate}
\begin{eqnarray}
&&\mathbb{E}\int_0^T (T-s)^{-2\beta}\|\mathcal{Y}^\gamma_s\|^2_H\, ds\nonumber\\
&\le& C \mathbb{E}\int_0^T (T-t)^{-2\beta}\int_0^t(t-s)^{-2\alpha}(T-s)^{2\alpha}\|\gamma_s\|^2_H\, ds\, dt\nonumber\\
&=&C \mathbb{E}\int_0^T \int_s^T(T-t)^{-2\beta}(t-s)^{-2\alpha}\,dt\, (T-s)^{2\alpha}\|\gamma_s\|^2_H\, ds\nonumber\\
&=&C \mathbb{E}\int_0^T \int_0^1 (1-\theta)^{-2\beta}\theta^{-2\alpha}\,d\theta\, (T-s)^{1-2\beta}\|\gamma_s\|^2_H\, ds\nonumber\\
&\le& C T^{1-2\beta}\int_0^1 (1-\theta)^{-2\beta}\theta^{-2\alpha}\,d\theta \, \|\gamma\|^2_{L^2_{\cal P}(\Omega\times [0,T],H)}.
\end{eqnarray}
Here we have used in the  last equality the transformation of variables: $t=s+(T-s)\theta$.
\qed

\begin{remark} \label{R 2.3} Let $N\ge 1$ be an integer.  For $\gamma\in L^\infty_{\cal P}(\Omega\times [0,T],H)$, the linear stochastic equation
\begin{equation}
    \left\{\begin{array}{lll}
    d\mathcal{Y}_t&=&\displaystyle A\mathcal{Y}_t\,dt + \sum_{i=1}^\infty C_i(t)\mathcal{Y}_t\,d\beta_t^i+\sum_{i=1}^N C_i(t)(T-t)^\alpha\gamma_t\,d\beta_t^i,\quad t\in (0,T];
    \\
    \mathcal{Y}_0&=&0,
\end{array}\right.
\end{equation}
has a unique mild solution $\mathcal{Y}^{\gamma, N}$ in $L^2_{\cal P}(\Omega, C([0,T],H)).$ Furthermore, we have
\begin{equation}\label{max estimate-N}
    \mathbb{E}\|\mathcal{Y}^{\gamma,N}_t\|^2_H\le C  \mathbb{E}\int_0^t(t-s)^{-2\alpha}(T-s)^{2\alpha}\|\gamma_s\|^2_H\, ds
\end{equation}
for a positive constant $C$, which does not depend on $N$.
For $\eta\in L^2(\Omega, \mathcal{F}_T,\mathbb{P}, H)$ and $f_0\in L^0_{\cal P}(\Omega\times [0,T],H)$ such that for some $\beta\in (0,\frac{1}{2})$
\begin{equation}\label{f0-2}
   {\mathbb E}\int_0^T(T-s)^{2\beta}\|f_0(s)\|^2_H\, ds<\infty,
\end{equation}
the following linear functional $G^N$ defined by
\begin{equation}\label{LF-2}
  G^N(\gamma):=E\langle \eta, \mathcal{Y}^{\gamma, N}_T\rangle+E\int_0^T\langle f_0(t), \mathcal{Y}^{\gamma,N}_t\rangle\, dt, \quad \gamma\in L^\infty_{\cal P}(\Omega\times [0,T],H)
\end{equation}
has a unique  linear and continuous extension to $L^2_{\cal P}(\Omega\times [0,T],H)$, denoted by  $\overline G^N$. Furthermore, there is a positive constant $C$ such that $C$ does not depend on $N$ and
\begin{eqnarray}
|\overline G^N(\gamma)|^2&\le& C\|\gamma\|^2_{L^2_{\cal P}(\Omega\times [0,T],H)}, \quad \gamma\in L^2_{\cal P}(\Omega\times [0,T],H).
\end{eqnarray}
\end{remark}

We now recall the following result from \cite{FHT}, concerning the following linear BSEE driven by a white noise:
\begin{equation}\label{lbsde}
    \left\{\begin{array}{lll}
    -dP_t&=&\displaystyle [A^*P_t + \sum_{i=1}^\infty C_i^*(t)Q_te_i+f_0(t)]\, dt-\sum_{i=1}^\infty Q_te_i\,d\beta_t^i,
    \\
    P_T&=&\eta.
\end{array}\right.
\end{equation}

\begin{lemma} \label{2.1} Let $\eta\in L^2(\Omega,{\cal F}_T,\mathbb P,H)$ and $f_0\in L^2_{\cal P}(\Omega\times [0,T],H)$.
There exists a unique mild solution $(P,Q)\in  L^2_{{\cal P}}(\Omega\times [0,T],H)\times L^2_{{\cal P}}(\Omega\times [0,T],\mathcal{S}_2(H))$ to BSEE~(\ref{lbsde}).
\end{lemma}

In the subsequent study of the nonlinear case, we need the following a priori estimate for BSEE~(\ref{lbsde}).

\begin{proposition}(a priori estimate) Let $\eta\in L^2(\Omega,{\cal F}_T,\mathbb P,H)$ and $f_0\in L^2_{\cal P}(\Omega\times [0,T],H)$. For $\beta\in (\alpha,\frac{1}{2})$ and a mild solution $(P,Q)\in  L^2_{{\cal P}}(\Omega\times [0,T],H\times\mathcal{S}_2(H))$ to BSEE~(\ref{lbsde}), we have
\label{apriori}
\begin{eqnarray*}
& &\mathbb E\int_t^T ||P_s||_H^2 ds+\mathbb E\int_t^T||Q_s||^2_{{\mathcal S}_2(H)}ds+\mathbb E\int_t^T (T-s)^{2\alpha}||\sum_{i=1}^\infty C_i(s)Q_se_i||_H^2 ds\\
&\le&C\left( \mathbb E||\eta||_H^2 +\mathbb E\int_t^T(T-s)^{2\beta} || f_0||_H^2 ds\right).
\end{eqnarray*}
\end{proposition}

\noindent{\bf Proof.} Let us prove the first two terms by duality argument.
We have,
\begin{eqnarray*}
\mathbb E\int_t^T\langle P_s,\rho_s\rangle_H ds+\mathbb E\int_t^T\langle Q_s,\Gamma_s\rangle_{\mathcal{S}_2(H)}\, ds
&=&\mathbb E \langle\eta,X_T\rangle_H +\mathbb E\int_t^T\langle f_s,X_s\rangle_H ds,
\end{eqnarray*}
where
$$dX_s=(AX_s+\rho_s)ds+\sum_{i=1}^\infty C_i(s)X_sd\beta_s^i+\sum_{i=1}^\infty\Gamma_se_id\beta_s^i, \quad s\in (t,T]; \quad X_t=0.$$
Hence,
\begin{eqnarray*}
& &\mathbb E\int_t^T\langle P_s,\rho_s\rangle_H ds+\mathbb E\int_t^T\langle Q_s,\Gamma_s\rangle_{\mathcal{S}_2(H)}ds\\
&=&\mathbb E \langle\eta,X_T\rangle_H +\mathbb E\int_t^T\langle f_s,X_s\rangle_H ds\\
&\le&\mathbb E\left[||\eta||^2_H\right]^{\frac{1}{2}}\mathbb E[||X_T||^2_H]^{\frac{1}{2}}
+\int_t^T {\mathbb E}[||X_s||^2_H]^{\frac{1}{2}}\mathbb{E}[||f_s||^2_H]^{\frac{1}{2}}\, ds\\
&\le&\mathbb E\left[||\eta||^2_H\right]^{\frac{1}{2}}\mathbb E[||X_T||^2_H]^{\frac{1}{2}}
+\sup_{t\le s\le T}{\mathbb E}[||X_s||^2_H]^{\frac{1}{2}}\int_t^T \mathbb{E}[||f_s||^2_H]^{\frac{1}{2}}\, ds\\
&\le&\mathbb E\left[||\eta||^2_H\right]^{\frac{1}{2}}\mathbb E[||X_T||^2_H]^{\frac{1}{2}}
+\sup_{t\le s\le T}{\mathbb E}[||X_s||^2_H]^{\frac{1}{2}}\int_t^T (T-s)^\beta\mathbb{E}[||f_s||^2_H]^{\frac{1}{2}}(T-s)^{-\beta}\, ds\\
&\le&\mathbb E\left[||\eta||^2_H\right]^{\frac{1}{2}}\mathbb E[||X_T||^2_H]^{\frac{1}{2}}
+{\frac{T^{1-2\beta}}{1-2\beta}}\sup_{t\le s\le T}{\mathbb E}[||X_s||^2_H]^{\frac{1}{2}}\mathbb{E}\left[\int_t^T (T-s)^{2\beta} ||f_s||^2_H\, ds\right]^{\frac{1}{2}}.
\end{eqnarray*}
Using the inequality from \cite[Proposition~4.4]{FHT}
$$\sup_{s\in [t,T]}\mathbb E\|X_s\|^2_H\le C\|(\rho, \Gamma)\|^2_{L^2_{{\cal P}}(\Omega\times [t,T],H\times \mathcal{S}_2(H))}$$
with
$$
\|(\rho, \Gamma)\|^2_{L^2_{{\cal P}}(\Omega\times [t,T],H\times \mathcal{S}_2(H))}:=\mathbb E\int_t^T\|\rho_s\|^2_Hds+\mathbb E\int_t^T\|\Gamma_s\|_{{\mathcal S}_2(H)}^2ds,
$$
we have
\begin{eqnarray*}
& &\mathbb E\int_t^T\langle P_s,\rho_s\rangle_H ds+\mathbb E\int_t^T\langle Q_s,\Gamma_s\rangle_{{\mathcal S}_2(H)}ds\\
&\le&C\left(\mathbb E\|\eta\|^2_H
+\mathbb{E}\int_t^T (T-s)^{2\beta} ||f_s||^2_H\, ds\right)^{\frac{1}{2}}\|(\rho, \Gamma)\|_{L^2_{{\cal P}}(\Omega\times [t,T],H\times \mathcal{S}_2(H))}.
\end{eqnarray*}
This implies the desired a priori estimate for the norm of $(P, Q)$.

To complete the proof, we consider once again the duality:
\begin{eqnarray*}
\mathbb E\int_t^T\langle(T-s)^\alpha\sum_{i=1}^\infty  C_iQ_se_i,\gamma_s\rangle ds
&=&\mathbb E \langle\eta,X_T\rangle +\mathbb E\int_t^T\langle f_s,X_s\rangle ds,
\end{eqnarray*}
where $X_t=0$ and
$$dX_s=AX_sds+\sum_{i=1}^\infty C_i(s)X_sd\beta_s^i+\sum_{i=1}^\infty C_i(s)(T-s)^\alpha\gamma_sd\beta_s^i, \quad s\in [t,T].$$
We have
\begin{eqnarray*}
& &\mathbb E\int_t^T\langle(T-s)^\alpha\sum_{i=1}^\infty  C_iQ_se_i,\gamma_s\rangle ds\\
&=&\mathbb E \langle\eta,X_T\rangle +\mathbb E\int_t^T\langle(T-s)^{\beta} f_s, (T-s)^{-\beta} X_s\rangle ds\\
&\le&\mathbb E\left[||\eta||^2_H\right]^{\frac{1}{2}}\mathbb E\left[||X_T||^2_H\right]^{\frac{1}{2}} \\
&&+\mathbb E\left[\int_t^T(T-s)^{2\beta}|| f_s||^2_Hds\right]^{\frac{1}{2}}\mathbb E\left[\int_t^T(T-s)^{-2\beta} ||X_s||^2_H ds\right]^{\frac{1}{2}}.
\end{eqnarray*}

On the one hand, we have from \cite[Proposition 4.5]{FHT},
\begin{eqnarray*}
\mathbb E\|X_T\|^2_H&\le& C\int_t^T(T-l)^{-2\alpha}(T-l)^{\alpha}\mathbb E\|\gamma_l\|_H^2dl= C\int_t^T\mathbb E\|\gamma_l\|_H^2dl
\end{eqnarray*}
and
\begin{eqnarray*}
& &\mathbb E\int_t^T(T-s)^{-2\beta} ||X_s||^2_H ds\\
&\le&C\int_t^T(T-s)^{-2\beta} \int_t^s(s-l)^{-2\alpha}(T-l)^{2\alpha}\mathbb E\|\gamma_l\|_H^2dlds\\
&=&C\int_t^T \int_l^T (T-s)^{-2\beta} (s-l)^{-2\alpha}ds(T-l)^{2\alpha}\mathbb E\|\gamma_l\|_H^2dl\\
&=&C\left(\int_0^1(1-\theta)^{-2\beta}\theta^{-2\alpha}d\theta\right)\int_t^T(T-l)^{1-2\beta}\mathbb E\|\gamma_l\|_H^2dl\\
&\le&CT^{1-2\beta}\left(\int_0^1(1-\theta)^{-2\beta}\theta^{-2\alpha}d\theta\right)\int_t^T\mathbb E\|\gamma_l\|_H^2dl.
\end{eqnarray*}
Note that in the last equality, we have used the transformation of variables: $s=l+(T-l)\theta$.
Concluding the above, we have
\begin{eqnarray*}
&&\mathbb E\int_t^T\langle(T-s)^\alpha\sum_{i=1}^\infty  C_iQ_se_i,\gamma_s\rangle ds\\
&\le &C \left(\mathbb E||\eta||^2_H
+\mathbb{E}\int_t^T (T-s)^{2\beta} ||f_s||^2_H\, ds\right)^{\frac{1}{2}}\|\gamma\|_{L^2_{{\cal P}}(\Omega\times [t,T],H)}.
\end{eqnarray*}
Then we have the last desired a priori estimate. \qed

\begin{proposition}\label{2.3} Suppose that
$$ \mathbb E||\eta||_H^2 +\mathbb E\int_0^T(T-s)^{2\beta} || f_0||_H^2 ds<\infty.$$
There exists a unique solution to linear BSEE~\eqref{lbsde} such that
$$\mathbb E\int_0^T ||P_s||_H^2 ds+\mathbb E\int_0^T||Q_s||^2_{{\mathcal S}_2(H)}ds+\mathbb E\int_0^T (T-s)^{2\alpha}||\sum_{i=1}^\infty C_i(s)Q_se_i||_H^2 ds<\infty.$$
\end{proposition}
\noindent{\bf Proof.} Uniqueness is an immediate consequence of Proposition~\ref{apriori}. It remains to consider the existence assertion.

Define for any integer $k>T$,
$$
f_0^k(t)= f_0(t)\chi_{[0, T-1/k]}(t), \quad t\in [0,T].
$$
We have for any $k>T$,
\begin{eqnarray}
  &&k^{-2\beta}E\int_0^T\|f_0^k(s)\|^2_H\, ds=k^{-2\beta}E\int_0^{T-1/k}\|f_0(s)\|^2_H\, ds\nonumber\\
  \le&& E\int_0^{T-1/k}(T-s)^{2\beta}\|f_0(s)\|^2_H\, ds\le  E\int_0^T (T-s)^{2\beta}\|f_0(s)\|^2_H\, ds<\infty.
\end{eqnarray}
Therefore, $f_0^k\in L^2_{\cal P}(\Omega\times [0,T],H)$, and in view of Lemma~\ref{2.1},  BSEE~\eqref{lbsde} for $f_0=f_0^k$ has a unique mild solution $(P^k, Q^k)$ for any integer $k>T$.

Moreover, for $k>l>T$, we have
\begin{eqnarray}
  E\int_0^T(T-s)^{2\beta}\|f_0^k(s)-f_0^l(s)\|^2_H\, ds \le E\int_{T-1/l}^{T-1/k}(T-s)^{2\beta}\|f_0(s)\|^2_H\, ds\rightarrow 0
\end{eqnarray}
as $k,l\rightarrow \infty$. Applying Proposition~\ref{apriori}, we see that $\{(P^k,Q^k)\}$ is a Cauchy sequence in the space $L^2_{\cal P}(\Omega\times [0,T],H\times \mathcal{S}_2(H))$, and the sequence of processes  $\{(T-s)^{\alpha}\sum_{i=1}^\infty C_i(s)Q^k_se_i, s\in [0,T]; k>T\} $ is a Cauchy sequence in the space $L^2_{\cal P}(\Omega\times [0,T],H)$. Thus, they have limits $(P,Q,(T-t)^\alpha S)\in L^2_{\cal P}(\Omega\times [0,T],H\times \mathcal{S}_2(H)\times H)$, which satisfies the following equation:
 \begin{equation}
    \left\{\begin{array}{lll}
    -dP_t&=&\displaystyle [A^*P_t + S(t)+f_0(t)]\, dt-\sum_{i=1}^\infty Q_te_i\,d\beta_t^i,
    \\
    P_T&=&\eta.
\end{array}\right.
\end{equation}

To show that $(P,Q)$ is a mild solution to BSEE~\eqref{lbsde}, it is sufficient for us to prove that
\begin{equation}\label{desired}
  S(t)=\sum_{i=1}^\infty C_i^*(t)Q_te_i,
\end{equation}
with the limit being defined in the following weak sense:
$$
\{(T-t)^\alpha\sum_{i=1}^N C_i^*(t)Q_te_i, t\in [0,T]\} \mathop{\longrightarrow}^{N\to \infty}  \{(T-t)^\alpha S(t), t\in [0,T]\}
$$
weakly in the Hilbert space $L^2_{\cal P}(\Omega\times [0,T],H)$.

Note that $\mathcal{Y}^\gamma$ is the solution to the stochastic equation~\eqref{lfsde} for $\gamma\in L^\infty_{\cal P}(\Omega\times [0,T],H)$. We have the following duality:
\begin{eqnarray}
  &&E\int_0^T\langle (T-t)^\alpha\sum_{i=1}^\infty C_i^*(t)Q^ke_i , \gamma_t \rangle\, dt =  E\int_0^T\sum_{i=1}^\infty\langle Q^ke_i , (T-t)^\alpha C_i(t)\gamma_t \rangle\, dt\nonumber\\
  &=& E\langle \eta, \mathcal{Y}^\gamma_T\rangle+E\int_0^T\langle f_0^k(t), \mathcal{Y}^\gamma_t\rangle\, dt.
\end{eqnarray}
Passing to the limit $k\to \infty$, we have for $\gamma\in L^\infty_{\cal P}(\Omega\times [0,T],H)$,
\begin{eqnarray}\label{S}
  E\int_0^T\langle (T-t)^\alpha S(t), \gamma_t \rangle\, dt
  &=& E\langle \eta, \mathcal{Y}^\gamma_T\rangle+E\int_0^T\langle f_0(t), \mathcal{Y}^\gamma_t\rangle\, dt=G(\gamma).
\end{eqnarray}
Since the process $(T-t)^\alpha S(t), t\in [0,T]$ lies in $L^2_{\cal P}(\Omega\times [0,T],H)$, in view of Lemma~\ref{2.2}, we have for $\gamma\in L^2_{\cal P}(\Omega\times [0,T],H)$
\begin{eqnarray}\label{S-ex}
  E\int_0^T\langle (T-t)^\alpha S(t), \gamma_t \rangle\, dt
  &=& {\overline G}(\gamma).
\end{eqnarray}

On the other hand,  we have the duality:
\begin{eqnarray}
  &&E\int_0^T\left\langle (T-t)^\alpha\sum_{i=1}^N C_i^*(t)Q^ke_i , \gamma_t \right\rangle\, dt =  E\int_0^T\sum_{i=1}^N\langle Q^ke_i , (T-t)^\alpha C_i(t)\gamma_t \rangle\, dt\nonumber\\
  &=& E\left\langle \eta, \mathcal{Y}^N_T\right\rangle+E\int_0^T\left\langle (T-t)^\beta f_0^k(t), (T-t)^{-\beta}\mathcal{Y}^N_t\right\rangle\, dt.
\end{eqnarray}
Setting $k\to \infty$, in view of \eqref{max estimate-N}, we have
\begin{eqnarray}\label{dual-N}
  E\int_0^T\left\langle (T-t)^\alpha S^N(t), \gamma_t \right\rangle\, dt
  &=& E\left\langle \eta, \mathcal{Y}^N_T\right\rangle+E\int_0^T\left\langle f_0(t), \mathcal{Y}^N_t\right\rangle\, dt=G^N(\gamma)
\end{eqnarray}
with
\begin{equation}\label{SN}
    S^N(t):=\sum_{i=1}^N C_i^*(t)Q e_i, \quad t\in [0,T].
\end{equation}
In view of Remark \ref{R 2.3}, we have
\begin{equation}\label{}
 \mathbb{E}\int_0^T \|(T-t)^\alpha S^N(t)\|^2_H\, dt= \|\overline G^N(\cdot)\|^2\le C
\end{equation}
with $C$ being independent of $N$. Then, the set $\{(T-\cdot)^\alpha S^N(\cdot), N=1,2,\ldots\}$ is  weakly compact, and thus has a weakly convergent subsequence.
Let $(T-\cdot)\overline S(\cdot)$ be one weak limit. Then in view of equality~\eqref{dual-N}, we have for $\gamma\in L^2_{\cal P}(\Omega\times [0,T],H)$,
\begin{eqnarray}
E\int_0^T\left\langle (T-t)^\alpha {\overline S}(t), \gamma_t \right\rangle\, dt
  &=& \overline G(\gamma)=E\int_0^T\left\langle (T-t)^\alpha S(t), \gamma_t \right\rangle\, dt.
\end{eqnarray}
Therefore, $\overline S=S$, and  the desired equality~\eqref{desired} is true. \qed

\section{Main Result}

In this section, we state and prove the following result.

\begin{theorem}\label{Main result} For $\eta\in L^2(\Omega, \mathcal{F}_T,\mathbb{P}, H)$ and the map $f: \Omega\times [0,T]\times H\times \mathcal{S}_2(H)\times H\to H$ where $f(\cdot,0,0,0)\in L^0_{{\cal P}}(\Omega\times [0,T],H)$ such that 
\begin{equation}\label{f(0,0,0)}
   {\mathbb E}\int_0^T(T-s)^{2\beta}\|f(s,0,0,0)\|^2_H\, ds<\infty
\end{equation}
for some $\beta\in (\alpha,\frac{1}{2})$, and 
\begin{equation}\label{Lip}
\|f(t,p_1,q_1,s_1)-f(t,p_2,q_2,s_2)\|_H\le L(\|p_1-p_2\|_H+\|q_1-q_2\|_{{\mathcal S}_2(H)}+\|s_1-s_2\|_H)
\end{equation}
for a positive constant $L$. 
There exists a unique solution $(P,Q)$ for (\ref{nlbsde})
such that
$$\mathbb E\int_0^T ||P_s||_H^2 ds+\mathbb E\int_0^T||Q_s||^2_{{\mathcal S}_2(H)}ds+\mathbb E\int_0^T (T-s)^{2\alpha}||\sum_{i=1}^\infty C_i(s)Q_se_i||_H^2 ds<\infty.$$
\end{theorem}

\noindent{\bf Proof.} (i) Uniqueness. Let $(P^i,Q^i)$ be a mild solution to BSEE~\eqref{nlbsde} for $i=1,2$. Define
$$
\widehat P:=P^1-P^2, \quad \widehat Q:=Q^1-Q^2;
$$
and
$$ {\widehat f}(t):=f(t,P^1_t,Q^1_t,\sum_{i=1}^\infty C_i^*(t)Q^1_te_i)-f(t,P^2_t,Q^2_t,\sum_{i=1}^\infty C_i^*(t)Q^2_te_i), \quad t\in [0,T].
$$ We have
\begin{equation}\label{lbsde1}
    \left\{\begin{array}{lll}
    -d{\widehat P}_t&=&\displaystyle [A^*{\widehat P}_t + \sum_{i=1}^\infty C_i^*(t){\widehat Q}_te_i+{\widehat f}(t)]\,dt-\sum_{i=1}^\infty {\widehat Q}_te_i\,d\beta_t^i, \quad t\in [0,T];
    \\
   {\widehat P}_T&=&0.
\end{array}\right.
\end{equation}
It suffices to show that $\widehat P=0$ and $\widehat Q=0$ on the interval $[T-\varepsilon_0, T]$ for a sufficiently small $\varepsilon_0>0$.  From Proposition~\ref{apriori}, we have
\begin{eqnarray*}
& &\mathbb E\int_t^T ||{\widehat P}_s||_H^2 ds+\mathbb E\int_t^T||{\widehat Q}_s||^2_{{\mathcal S}_2(H)}ds+\mathbb E\int_t^T (T-s)^{2\alpha}||\sum_{i=1}^\infty C_i(s){\widehat Q}_se_i||_H^2 ds\\
&\le&C\left( \mathbb E\int_t^T(T-s)^{2\beta} || {\widehat f}(s)||_H^2 ds\right),
\end{eqnarray*}
and further in view of the Lipschitz continuity of $f$,
\begin{eqnarray*}
& &\mathbb E\int_t^T ||{\widehat P}_s||_H^2 ds+\mathbb E\int_t^T||{\widehat Q}_s||^2_{{\mathcal S}_2(H)}ds+\mathbb E\int_t^T (T-s)^{2\alpha}||\sum_{i=1}^\infty C_i(s){\widehat Q}_se_i||_H^2 ds\\
&\le& C \mathbb E\int_t^T(T-s)^{2(\beta-\alpha)} \left(\| {\widehat P}_s\|_H^2+ ||{\widehat Q}_s||^2_{{\mathcal S}_2(H)}+||(T-s)^{2\alpha}\sum_{i=1}^\infty C_i(s){\widehat Q}_te_i||_H^2 \right)\, ds\\
&\le& C \varepsilon^{2(\beta-\alpha)}\mathbb E\int_t^T \left(\| {\widehat P}_s\|_H^2+ ||{\widehat Q}_s||^2_{{\mathcal S}_2(H)}+||(T-s)^{2\alpha}\sum_{i=1}^\infty C_i(s){\widehat Q}_se_i||_H^2 \right)\, ds.
\end{eqnarray*}
Thus we have the desired uniqueness on the interval $[T-\varepsilon_0, T]$ for a sufficiently small $\varepsilon_0>0$. Iteratively in a backward way, we can show the uniqueness on the whole interval $[0,T]$.

We use the Picard iteration to construct a sequence of solutions to linear BSEEs, and show that its limit is a solution to the nonlinear BSEE~\eqref{nlbsde}.
Noting that $f(\cdot,0,0,0)$ verifies the integrability~\eqref{f(0,0,0)},  in view of Proposition~\ref{2.3}, the following BSEE
\begin{equation}\label{0bsde1}
    \left\{\begin{array}{lll}
    -dP_t^1&=&\displaystyle [A^*P^1_t + \sum_{i=1}^\infty C_i^*(t)Q^1_te_i]\,dt+f(t,0,0,0)\,dt  -\sum_{i=1}^\infty Q^1_te_i\,d\beta_t^i, \quad t\in [0,T];
    \\
    P^1_T&=&\eta
\end{array}\right.
\end{equation}
has a unique mild solution $(P^1, Q^1)$, and BSEE
\begin{equation}\label{aprox bsde1}
    \left\{\begin{array}{lll}
    -dP_t^{k+1}&=&\displaystyle [A^*P^{k+1}_t + \sum_{i=1}^\infty C_i^*(t)Q^{k+1}_te_i]\,dt+f(t,P^k_t,Q^k_t,\sum_{i=1}^\infty C_i^*(t)Q^k_te_i)\ dt\\
    &&\displaystyle -\sum_{i=1}^\infty Q^{k+1}_te_i\,d\beta_t^i, \quad t\in [0,T];
    \\
    P^{k+1}_T&=&\eta,
\end{array}\right.
\end{equation}
has a unique mild solution $(P^{k+1}, Q^{k+1})$ with $k=1,2,\ldots,$   such that for  $k=0,1,2,\ldots,$
$$\mathbb E\int_0^T ||P^{k+1}_s||_H^2 ds+\mathbb E\int_0^T||Q^{k+1}_s||^2_{{\mathcal S}_2(H)}ds+\mathbb E\int_0^T (T-s)^{2\alpha}||\sum_{i=1}^\infty C_i(s)Q^{k+1}_se_i||_H^2 ds<\infty.$$

From Proposition~\ref{apriori} and the Lipschitz continuity of $f$, we can show the following for $t\in [T-\varepsilon, T]$
\begin{eqnarray*}
&&\mathbb E\int_t^T ||P^{k+1}_s-P^k_s||_H^2 ds+\mathbb E\int_t^T||Q^{k+1}_s-Q^k_s||^2_{{\mathcal S}_2(H)}ds\\
&&\quad+\mathbb E\int_t^T (T-s)^{2\alpha}||\sum_{i=1}^\infty C_i(s)[Q^{k+1}_s-Q^k_s]e_i||_H^2 ds\\
&\le &C \varepsilon^{2(\beta-\alpha)}\left(\mathbb{E}\int_t^T ||P^{k-1}_s-P^k_s||_H^2 ds+\mathbb E\int_t^T||Q^{k-1}_s-Q^k_s||^2_{{\mathcal S}_2(H)}ds\right)\\
&&\quad +C \varepsilon^{2(\beta-\alpha)}\mathbb E\int_t^T (T-s)^{2\alpha}||\sum_{i=1}^\infty C_i(s)[Q^{k-1}_s-Q^k_s]e_i||_H^2\,ds.
\end{eqnarray*}
Choose $\varepsilon_1>0$ such that $C \varepsilon_1^{2(\beta-\alpha)}={\frac{1}{2}}$. Then the sequence
$$(P^k_s,Q^k_s, (T-s)^{2\alpha}||\sum_{i=1}^\infty C_i(s)Q^k_se_i||_H^2), \quad s\in [T-\varepsilon_1, T]$$
converges  strongly in $L^2_{{\cal P}}(\Omega\times [T-\varepsilon_1, T],H\times \mathcal{S}_2(H)\times H)$ to a triplet
$$\{(P_t, Q_t, (T-t)^\alpha S(t)), t\in [T-\varepsilon_1,T]\}.$$
Moreover, we have
\begin{equation}\label{limit bsde1}
    \left\{\begin{array}{lll}
    -dP_t&=&\displaystyle [A^*P_t + S(t)+f(t,P_t,Q_t,S(t))]\,dt -\sum_{i=1}^\infty Q_te_i\,d\beta_t^i, \quad t\in [0,T];    \\
    P_T&=&\eta.
\end{array}\right.
\end{equation}
It remains to prove the following weak convergence:
\begin{equation}\label{weak limit}
\{(T-t)^\alpha\sum_{i=1}^N C_i^*(t)Q_te_i, t\in [T-\varepsilon_1,T]\} \mathop{\longrightarrow}^{N\to \infty}  \{(T-t)^\alpha S(t), t\in [T-\varepsilon_1,T]\}
\end{equation}
weakly in the Hilbert space $L^2_{\cal P}(\Omega\times [T-\varepsilon_1,T],H)$.

For $\gamma\in L^\infty_{\cal P}(\Omega\times [T-\varepsilon_1,T],H),$ let $\mathcal{X}^\gamma$ be the unique mild solution of the following stochastic equation:
\begin{equation}
    \left\{\begin{array}{lll}
    d\mathcal{X}_t&=&\displaystyle A\mathcal{X}_t\,dt + \sum_{i=1}^\infty C_i(t)\mathcal{X}_t\,d\beta_t^i+\sum_{i=1}^\infty C_i(t)(T-t)^\alpha\gamma_t\,d\beta_t^i,\quad t\in (T-\varepsilon_1,T];
    \\
    \mathcal{X}_{T-\varepsilon_1}&=&0.
\end{array}\right.
\end{equation}
Using~\eqref{aprox bsde1}, we have the following duality:
\begin{eqnarray}
 && E\int_{T-\varepsilon_1}^T\langle (T-t)^\alpha  \sum_{i=1}^\infty C_i^*(t)Q^{k+1}_te_i, \gamma_t \rangle\, dt\nonumber\\
  &=& E\langle \eta, \mathcal{X}^\gamma_T\rangle+E\int_{T-\varepsilon_1}^T\langle f(t, P^k_t,Q^k_t, \sum_{i=1}^\infty C_i^*(t)Q^k_te_i), \mathcal{X}^\gamma_t\rangle\, dt.
\end{eqnarray}
By setting $k\to \infty$, we have
\begin{eqnarray}\label{nS}
  E\int_{T-\varepsilon_1}^T\langle (T-t)^\alpha S(t), \gamma_t \rangle\, dt
  &=& E\langle \eta, \mathcal{X}^\gamma_T\rangle+E\int_{T-\varepsilon_1}^T\langle f(t, P_t,Q_t,S(t)), \mathcal{X}^\gamma_t\rangle\, dt.
\end{eqnarray}

On the other hand, in view of  BSEE~\eqref{aprox bsde1},  we have the duality:
\begin{eqnarray}
  &&E\int_{T-\varepsilon_1}^T\left\langle (T-t)^\alpha\sum_{i=1}^N C_i^*(t)Q^{k+1}e_i , \gamma_t \right\rangle\, dt =  E\int_{T-\varepsilon_1}^T\sum_{i=1}^N\langle Q^{k+1}e_i , (T-t)^\alpha C_i(t)\gamma_t \rangle\, dt\nonumber\\
  &=& E\left\langle \eta, \mathcal{X}^{\gamma,N}_T\right\rangle+E\int_{T-\varepsilon_1}^T\left\langle f(t, P^k_t,Q^k_t,\sum_{i=1}^\infty C^*_i(t)Q^ke_i), \mathcal{X}^{\gamma,N}_t\right\rangle\, dt.
\end{eqnarray}
Here, $\mathcal{X}^{\gamma,N}$ is the unique mild solution of the following stochastic equation:
\begin{equation}
    \left\{\begin{array}{lll}
    d\mathcal{X}_t&=&\displaystyle A\mathcal{X}_t\,dt + \sum_{i=1}^\infty C_i(t)\mathcal{X}_t\,d\beta_t^i+\sum_{i=1}^N C_i(t)(T-t)^\alpha\gamma_t\,d\beta_t^i,\quad t\in (T-\varepsilon_1,T];
    \\
    \mathcal{X}_{T-\varepsilon_1}&=&0.
\end{array}\right.
\end{equation}
Setting $k\to \infty,$ we have
\begin{eqnarray}
  E\int_{T-\varepsilon_1}^T\left\langle (T-t)^\alpha\sum_{i=1}^N C_i^*(t)Qe_i , \gamma_t \right\rangle\, dt\nonumber
  &=& E\left\langle \eta, \mathcal{X}^{\gamma, N}_T\right\rangle+E\int_{T-\varepsilon_1}^T\left\langle f(t, P_t,Q_t,S(t)), \mathcal{X}^{\gamma,N}_t\right\rangle\, dt.
\end{eqnarray}
Note that we only have the following weaker integrability on $f(\cdot, P,Q,S)$:
$$\mathbb{E}\int_{T-\varepsilon_1}^T(T-t)^{2\alpha}\|f(t, P_t,Q_t,S(t))\|^2_H\, dt<\infty.$$
Subsequently, in view of \cite[Theorem 4.3]{FHT}, we have
\begin{eqnarray}
 &&\lim_{N\to \infty} E\int_{T-\varepsilon_1}^T\langle (T-t)^\alpha\sum_{i=1}^N C_i^*(t)Qe_i , \gamma_t \rangle\, dt\nonumber\\
  &=& E\langle \eta, \mathcal{X}^\gamma_T\rangle+E\int_{T-\varepsilon_1}^T\langle f(t, P_t,Q_t,S(t)), \mathcal{X}^\gamma_t\rangle\, dt.
\end{eqnarray}

Proceeding identically as in the proof of the previous equality~\eqref{desired}, we have the weak convergence~\eqref{weak limit}. In this way, we get the existence of BSEE~\eqref{nlbsde} on the interval $[T-\varepsilon_1, T]$. In a backward way, we can show its existence iteratively on the intervals $[T-2\varepsilon_1, T-\varepsilon_1], \ldots, [0, T-n_0\varepsilon_1]$ for the greatest integer $n_0$ such that $l\varepsilon_1<T$, and thus the existence on the whole interval $[0,T]$.
\qed

\section{Example}

Set $H:=L^2(0,1)$ and consider an  $H$-valued cylindrical Wiener process $\{ W_t, t\geq 0\}$. $A$ is the realization of the second derivative operator in $H$
with Dirichlet boundary conditions. So $\mathcal{D}(A)=H^2(0,1)\cap H^1_0(0,1)$ and $A\phi=\phi''$ for all $\phi \in \mathcal{D}(A)$. Choose an orthonormal basis in $L^2(0,1)$ with $\sup_i \sup_{x\in (0,1)}
|e_i(x)|<\infty$, for instance a trigonometrical  basis. Let $\sigma \in L^{\infty}_{\mathcal{P}}(\Omega\times (0,T), L^\infty(0,1))$. Define $C_i(t): H\to H$ by
$$
(C_i(t)\phi)(x):=\sigma(t,x)e_i(x)\phi(x), \quad (t, x)\in [0,T]\times [0,1]
$$
for $\phi\in H$.
We have $A^*=A$ and  $C_i^*=C_i$ with $i=1,2,\ldots$. From Da Prato and Zabczyk~\cite{DaZa2}, we see that $(A,C)$ satisfies Hypothesis~\ref{Hyp:C_i}.

Then for suitable conditions on $(\eta, f)$, our Theorem~\ref{Main result} can be applied to give the existence and uniqueness of a mild solution to the following backward stochastic partial differential equation driven by a space-time white noise:
\begin{equation}\label{Exam nlbsde}
    \left\{\begin{array}{lll}
    -dP_t(x)&=&\displaystyle \left[\frac{d^2} {dx^2}P_t(x) + \sigma(t,x)Q_t(x)\right]\, dt+f\left(t,x,P_t,Q_t,\sum_{i=1}^\infty \sigma(t)Q_t\right)\, dt\\
    &&\displaystyle -Q_t(x)\,dW_t(x), \quad (t,x)\in [0,T)\times (0,1);   \\
    P_t(0)&=&P_t(1)=0,\quad t\in [0,T];\\
    P_T(x)&=&\eta(x), \quad x\in [0,1].
\end{array}\right.
\end{equation}


\begin{thebibliography}{99}
\bibitem{Be}
A. Bensoussan.
 Stochastic maximum principle for distributed parameter systems.
J. Franklin Inst. 315 (1983), no. 5-6, 387--406.


\bibitem{DaZa} G. Da Prato, J.  Zabczyk.
Stochastic Equations in Infinite Dimensions.
Encyclopedia of Mathematics and its Applications, 44. Cambridge University Press, Cambridge,
 1992.

 \bibitem{DaZa2} G. Da Prato, J.  Zabczyk. Ergodicity for Infinite Dimensional Systems.
 London Mathematical Society Lecture Note Series, 229.
 Cambridge University Press, Cambridge, 1996.


\bibitem{DuMe2}
K. Du, Q. Meng.
A maximum principle for optimal control of stochastic evolution equations. SIAM J. Control Optim. 51
(2013), no. 6, 4343--4362.

 \bibitem{FuHuTe}
 M. Fuhrman, Y.  Hu, G. Tessitore.
 Stochastic maximum principle for optimal control of SPDEs.
 Appl. Math. Optim. 68 (2013), no. 2, 181--217.

\bibitem{FHT}
M. Fuhrman, Y. Hu  and G. Tessitore, Stochastic maximum principle for optimal control of partial differential equations driven by white noise.  arXiv:1409.4746v2,  2017.


\bibitem{Gu}
G. Guatteri.
 Stochastic maximum principle for SPDEs with noise and control on the boundary.
 Systems Control Lett. 60 (2011), no. 3, 198--204.



\bibitem{HuPe1}
Y. Hu, S. Peng.
Maximum principle for semilinear stochastic evolution control systems.
Stochastics Stochastics Rep. 33 (1990), no. 3-4, 159--180.

\bibitem{HuPe2} Y. Hu, S. Peng.
Adapted solution of a backward semilinear stochastic evolution equation.
Stochastic Anal. Appl. 9 (1991), no. 4, 445--459.




\bibitem{LuZh}
Q. L\"u, X. Zhang.
General Pontryagin-type stochastic maximum principle and
backward stochastic evolution equations in infinite dimensions.
SpringerBriefs in Mathematices. Springer, Cham, 2014.









\bibitem{TaLi}
S. Tang, X. Li.
 Maximum principle for optimal control of distributed parameter
 stochastic systems with random jumps.
 Differential equations, dynamical systems, and control science, 867--890,
 Lecture Notes in Pure and Appl. Math., 152, Dekker, New York,  1994.

%
%
%
%
%

\bibitem{Zh}
X.Y.  Zhou.
On the necessary conditions of optimal controls for stochastic
partial differential equations. SIAM J. Control Optim. 31 (1993), no. 6, 1462--1478.



\end{thebibliography}
\end{document}